\documentclass{article}
\usepackage{array}
\usepackage{amsmath}
\usepackage{amssymb}
\usepackage{graphicx}
\usepackage{listings}
\usepackage{url}

\usepackage[top=.9in,bottom=.9in]{geometry}
\newcommand{\sgn}{\mathop{\rm sgn}}
\newcommand{\csgn}{\mathop{\rm csgn}}
 \numberwithin{equation}{section}
 
\begin{document}
\pagestyle{plain}
\widowpenalty=600
\clubpenalty=600
\providecommand{\keywords}[1]
{
  \small	
  \textbf{\textit{Keywords---}} #1
}

\title{Application of the Argument Principle to Functions Expressed as Mellin Transforms}

\author{Bjoern S. Schmekel \footnote{bss28@cornell.edu} \\
\it Department of Physics, College of Studies for Foreign Diploma \\   
\it Recipients at the University of Hamburg, 20355 Hamburg, Germany }
\date{January 15, 2021}
\maketitle

\setcounter{page}{1}
\pagenumbering{arabic}

\begin{abstract}
We describe a numerical algorithm for evaluating the numbers of roots minus the number of poles contained in a region based on the argument principle
with the function of interest being written as a Mellin transformation of a usually simpler function. Because the function to be transformed may be simpler
than its Mellin transform whose roots are to be sought we express the final integrals in terms of the former accepting higher dimensional integrals. 
Nonlinear terms are expressed as convolutions approximating reciprocal values by exponential sums. As an example the final expression is applied to the Riemann Zeta function. 
The procedure is very inefficient numerically. However, depending on the function to be investigated it may be possible to find analytical estimates of the resulting integrals. 
\end{abstract}


\section{Introduction}
Object of this paper is to compute the number of roots minus the number of poles enclosed by a closed contour $C$ using the argument principle
\begin{eqnarray}
N_R - N_P = \frac{1}{2 \pi i} \oint_C \frac{f^\prime (s)}{f(s)} ds
\label{argprinciple}
\end{eqnarray}
where $f(z)$ is a meromorphic function on and inside of the contour $C$ which can be represented as an auxilliary function multiplied by the Mellin transform of another function 
\begin{eqnarray}
f(s) = K(s) Z(s) = K(s) \int_0^{\infty} z(t) t^{s-1} dt
\label{form}
\end{eqnarray}
assuming the Melling transform exists on and inside of the contour $C$. The latter is chosen such that it does not run over any poles or roots. 
We are interested in cases where $z(t)$ is a simple function therefore expressing the final result in terms of $K(s)$ as well as $z(t)$ instead of $Z(s)$. 
Furthermore, we are looking for an expression such that eqn. \ref{argprinciple} can be expressed as a multi-dimensional integral of $z(t)$. Obviously, the integrand in 
eqn. \ref{argprinciple} is a nonlinear functional of $z(t)$ (and $K(s)$ ), so arriving at such a result is not completely straightforward. In a first step we deal
with the nonlinearities by approximating the reciprocal value of $f(s)$ in the argument principle by an exponential sum \cite{Hackbusch2019,mclean2016exponential,McLean_2018,BEYLKIN2010131}, i.e.
\begin{eqnarray}
\frac{1}{x} \approx \sum_{j=1}^{N} \alpha_j e^{-c_j x} \equiv \mathcal{I}_1 (x)
\label{expsumapprox}
\end{eqnarray}
which is possible for $\Re(x)>0$. We will have to ensure this condition is always met possibly adding a factor which changes sign when appropriate. 

\begin{figure}
\scalebox{3.2}{\includegraphics[angle=0]{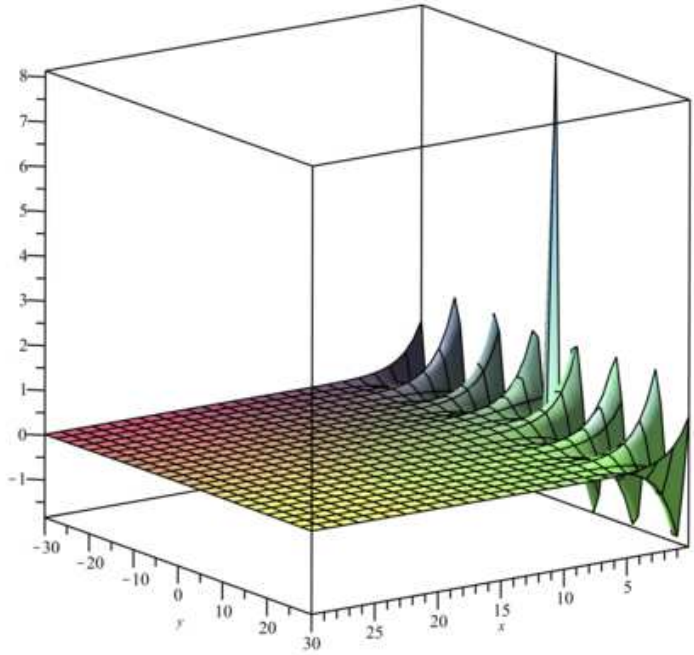}}
\scalebox{3.2}{\includegraphics[angle=0]{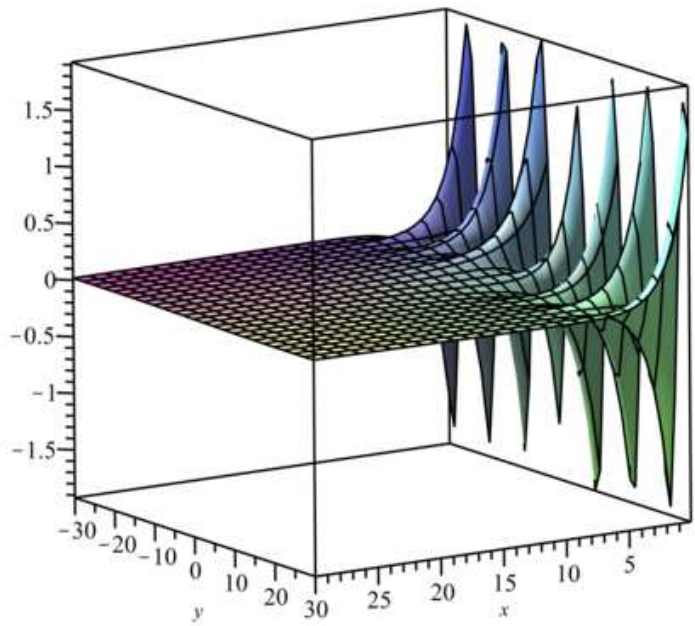}}
\caption{Real (above) and imaginary part (below) of the difference of $1/z$ and its exponential sum approximation with $z = x + iy$}
\label{accur1x}
\end{figure}

The exponential sum approximation has not been investigated too thoroughly for complex denominators. In fig. \ref{accur1x} The
approximation breaks down for small values of $\Re(x)$ as expected, but accuracy is not impacted by an imaginary part as only $\Re(x)>0$ is needed for convergence.

Using the following complex sign function
\begin{eqnarray}
\csgn(s) = 
\begin{cases}
-1 & \Re(s)<0 \\
1 & \Re(s)>0 \\
\sgn \left ( \Im(s) \right ) & \Re(s)=0
\end{cases}
\end{eqnarray}
and expanding the exponential function as a power series we obtain
\begin{eqnarray}
\nonumber
N_R - N_P & = & \frac{1}{2 \pi i} \oint_C ds \left [ K^{\prime} (s) Z(s) + K(s) Z^{\prime} (s) \right ] \sum_{j=1}^{N} \sum_{k=0}^{n} \alpha_j  \csgn \left (f(s) \right ) \frac{(-1)^k}{k!} c_j^k K^k(s) Z^k(s) {\csgn}^k \left (f(s) \right ) \\ \nonumber 
                   & = & \frac{1}{2 \pi i}  \sum_{j=1}^{N} \sum_{k=0}^{n} \oint_C ds   \alpha_j  {\csgn}^{k+1} \left (f(s) \right ) \frac{(-1)^k}{k!} c_j^k K^{\prime} (s) K^k(s) Z^{k+1}(s) \\ 
                   & + & \frac{1}{2 \pi i}  \sum_{j=1}^{N} \sum_{k=0}^{n} \oint_C ds   \alpha_j  {\csgn}^{k+1} \left (f(s) \right ) \frac{(-1)^k}{k!} c_j^k K^{k+1}(s) Z^{\prime} (s)  Z^k(s) 
                   = \int_0^{2\pi} d \phi \mathcal{K}(\phi)
\end{eqnarray}
The powers of $Z(s)$ can be expressed in terms of $z(t)$ using the Mellin convolution theorem
\begin{eqnarray}
Z(s) = \int_0^{\infty} dt  
z \left ( t \right )   t^{s-1} \\
Z^2(s) = \int_0^{\infty} dt  \int_0^{\infty} du_1 
z \left ( u_1 \right )  z \left ( \frac{t}{u_1} \right ) u_1^{-1}  t^{s-1} \\
Z^3(s) = \int_0^{\infty} dt  \int_0^{\infty} du_1  \int_0^{\infty}  du_2 
\label{Z3}
z \left ( u_1 \right ) z \left ( \frac{u_2}{u_1} \right )  z \left ( \frac{t}{u_2} \right ) u_1^{-1} u_2^{-1}  t^{s-1} \\
Z^k(s) = \int_0^{\infty} dt  \int_0^{\infty} du_1  \ldots \int_0^{\infty}  du_{k-1} t^{s-1}
z \left ( u_1 \right ) z \left ( \frac{t}{u_{k-1}} \right ) u_1^{-1} \prod_{j=1}^{k-2} z \left ( \frac{u_{j+1}}{u_j} \right )    u_{j+1}^{-1}
\label{Zk}
\end{eqnarray}
In appendix B a short Maple program is given which can be used to test the formulas given above. 
Similarly, exploiting standard rules for the Mellin transform
\begin{eqnarray}
Z^{\prime}(s) = \int_0^{\infty} dt  \ln(t) z \left ( t \right )   t^{s-1} \\
Z^{\prime} (s) Z(s) = \int_0^{\infty} dt  \int_0^{\infty} du_1 \ln (u_1)
z \left ( u_1 \right )  z \left ( \frac{t}{u_1} \right ) u_1^{-1}  t^{s-1} \\
\label{ZpZ}
\end{eqnarray}

If the behavior of the $\csgn$-function is non-trivial for the function to be investigated it may be approximated continuously by
\begin{eqnarray}
\csgn(x) \approx \tanh(x/\epsilon)
\label{csgntanh}
\end{eqnarray}
with precision increasing as $\epsilon \longrightarrow 0$ where $\tanh$ can be represented as
\begin{eqnarray}
\tanh(x) = - \frac{2i}{\pi} \int_0^{\infty} \frac{t^{\frac{2ix}{\pi}}-1}{t^2-1} dt
\label{intreptanh}
\end{eqnarray}
The integral converges if $-\pi/2 < \Im(x) < 0$.

\section{Example: Riemann Zeta Function}

We use the zeta function in the form \cite{milgram2012integral,Milgram_2013}
\begin{eqnarray}
\zeta(s) = \frac{2^{s-1}}{\left ( 1 - 2^{1-s} \right ) \Gamma(s+1)} \int_0^{\infty} \frac{t^s}{\cosh^2 (t) dt}
\label{defzeta}
\end{eqnarray}
which converges for $\Re(s)>-1$. The zeta function has been investigated using the argument principle before by many authors \cite{JOHNSON2009418}. 
The representation in eqn. \ref{defzeta} is written in the form of eqn. \ref{form}. Since the factor $K(s)$ in front of the integral has no roots or poles by itself
with the exception of the known pole at $s=1$ it would be sufficient to set $K(s)=1$. However, for the present purpose the full factor is retained in order to
stay in a numerically favorable range achieving sufficient accuracy in the exponential sum approximation. 

\begin{figure}
\scalebox{3.2}{\includegraphics[angle=270]{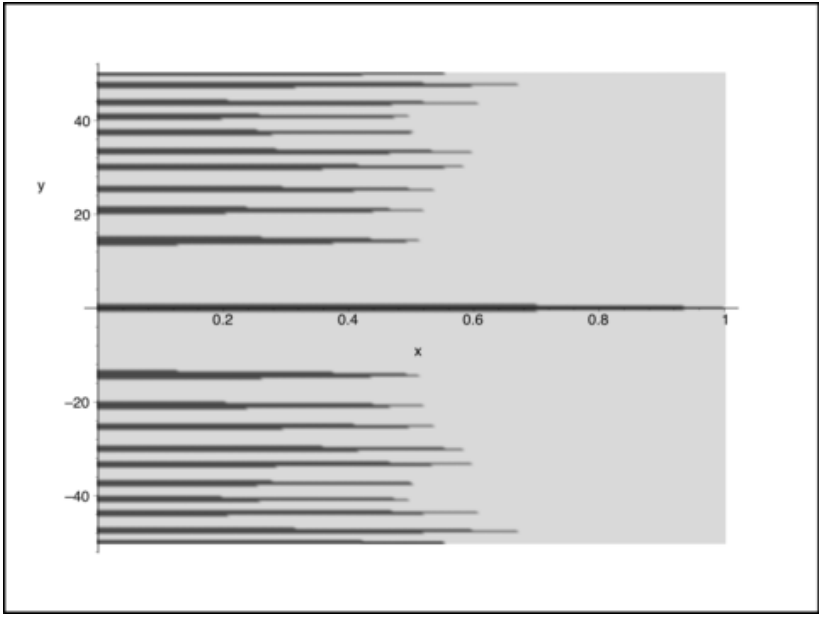}}
\caption{Regions of positive (white) and negative (black) real parts of the zeta function.}
\label{SgnReZeta}
\end{figure}

In table \ref{ResultsApproximations} we contrast the integrand in the argument principle with various steps towards the final approximation.
Integration is performed on a circle with radius $R=0.1$ around $z_0=0.57+1.57i$ not enclosing any roots. 
The second column is the integrand and factor of eqn. \ref{argprinciple} with the final $d \phi$-integration missing. The coefficients used in the exponential sum approximation
can be found in table \ref{coeffexpsumapprox} in the appendix. In the third column $1/f(z)$ is approximated by
$\mathcal{I}_1(z)$ which contains the exponential sum approximation. This is approximated further by $\mathcal{I}_2(z)$ where the exponential function has been
expanded in a power series up to linear order. Finally, in the fifth column the powers of $Z$ are expressed by Mellin convolutions given by eqn. \ref{Zk} and \ref{ZpZ}. 
Mathematica code producing the results in table  \ref{ResultsApproximations} can be found in appendix C. 

\begin{table}[]
\begin{tiny}
\begin{tabular}{rrrrr}
$\phi / (2 \pi)$ & $\frac{1}{2 \pi i} \frac{dz}{d \phi} \cdot \left .\frac{f'(z)}{f(z)} \right |_{z=z_0+R e^{i \phi}}$ & $\frac{dz}{d \phi} \cdot \left . \frac{f'(z) \mathcal{I}_1(z)}{2\pi i} \right |_{z=z_0+R e^{i \phi}}$ & $\frac{dz}{d \phi} \cdot \left . \frac{f'(z) \mathcal{I}_2(z)}{2 \pi i} \right |_{z=z_0+R e^{i \phi}}$ & $\mathcal{K}(\phi)$ \\
\hline
$0$    & $ 0.0124820  +0.0040853 i $  & $ 0.0155503   +0.0200828 i$  & $ 0.0155503 +0.0200828 i$  & $ 0.0155502+0.0200828 i$ \\
$1/8$ & $ 0.0062548  +0.0106734 i$ & $-0.0020676   +0.0219640 i$  & $-0.0020676 +0.0219640 i$  & $-0.0020676+0.0219640 i$ \\
$2/8$ & $-0.0021327  +0.0121535 i$ & $-0.0141169   +0.0148087 i$  & $-0.0141169 +0.0148087 i$  & $-0.0141169+0.0148087 i$ \\
$3/8$ & $-0.0101518  +0.0081128 i$ & $-0.0207191   +0.0033200 i$  & $-0.0207191 +0.0033200 i$ & $-0.0207191+0.0033200 i$ \\
$4/8$ & $-0.0140602  -0.0013970 i$ & $-0.0200646   -0.0121417 i$  & $-0.0200646 -0.0121417 i$  & $-0.0200647-0.0121417 i$ \\
$5/8$ & $-0.0089872  -0.0122548 i $ & $-0.0059828   -0.0263608 i$  & $-0.0059828 -0.0263608 i$ & $-0.0059828-0.0263609 i$ \\
$6/8$ & $ 0.0037589  -0.0148824 i$ & $ 0.0186236   -0.0229791 i$  & $ 0.0186236 -0.0229791 i$  & $ 0.0186236-0.0229791 i$ \\
$7/8$ & $ 0.0128362  -0.0064908 i$ & $ 0.0287771   +0.0013062 i$  & $ 0.0287771 +0.0013062 i$ & $ 0.0287771+0.0013062 i$ \\
$8/8$ & $ 0.0124820  +0.0040853 i $  & $ 0.0155503   +0.0200828 i$  & $ 0.0155503 +0.0200828 i$  & $  0.0155502+0.0200828 i$ 
\end{tabular}
\caption{$R=0.1$ and $z_0=0.57+1.57i$}
\label{ResultsApproximations}
\end{tiny}
\end{table}

Looking at the error introduced in each step we find that the expression by Mellin convolutions (cf. column 4 and 5 in table \ref{ResultsApproximations}) works fairly well.
The largest error is introduced by the exponential sum approximation which could be reduced by using more exponential terms (higher value of $N$ in eqn. \ref{expsumapprox}).
Ultimately, for arbitrarily high precision the number of terms in the expansion of the exponential function needs to be increased as well, though (higher value of $n$).  
Each new term introduces integrals of one more dimension which makes them increasingly hard to evaluate numerically. 

\section{Conclusions}
We presented a method which evaluates the number of roots minus the number of poles enclosed in a region using the argument principle focusing on function which
can be expressed as Mellin transforms of simple functions. The method was devised to work with the latter (simpler) function which was made possible by making use
of the exponential sum approximation and the expansion of the exponential function in a power series. The powers could be expressed in terms of Mellin convolutions of
the simpler function. 
Because of the high dimension of the involved integrals the method may not be feasible for high precision. However, since depending on the function of interest the integrands
may be simple it may be possible to come up with analytical estimates which may or may not exclude roots in a given region.

\section{Acknowledgments}
We acknowledge support by Wolfram Research having provided assistance with Mathematica and free maintenance thereof. 

\bibliography{bib}
\bibliographystyle{hunsrt}

\appendix
\section{Coefficients for the exponential sum approximation}

\begin{table}[h]
\begin{tabular}{rrr}
$i$ & $\alpha_i$ & $c_i$ \\
$1$ & $0.048$ & $0.017$ \\
$2$ & $0.235$ & $0.139$ \\
$3$ & $0.852$ & $0.627$ \\
$4$ & $2.737$ & $2.241$
\end{tabular}
\caption{Coefficients for the exponential sum approximation with values taken from \cite{Hackbusch2019}}
\label{coeffexpsumapprox}
\end{table}

\section{Maple Code Testing Eqn. \ref{Z3}} \footnote{Code tested using Maple 2019.2 for Mac OS X} 
The following Maple code computes the third power of the integral in eqn. \ref{defzeta} for $s=0.4$ without the factor in front using eqn. \ref{Z3} and by taking the third power directly. The results are
$0.4875296028$ and $0.4875296044$, respectively. For $s=0.4-0.3i$ we obtain $0.4103824778 + 0.1549090396i$ and $0.4103824766 + 0.1549090398i$, respectively. 

\begin{lstlisting}[language=C]
z:=unapply(t/cosh(t)^2,t);
integrand:=unapply(z(u2/u1)*z(t/u2)*z(u1)/u1/u2,u1,u2,t);
expr1:=Int(integrand(u1,u2,t),u1=0..infinity);
expr2:=Int(expr1,u2=0..infinity);
Int(expr2*t^(s-1),t=0..infinity);
subs(s=0.4,%);
evalf(%);
Zalt:=unapply(Zeta(s)*(1-2^(1-s))*GAMMA(s+1)/2^(s-1),s);
(Zalt(0.4))^3;
\end{lstlisting}

\section{Mathematica Code Producing Table \ref{ResultsApproximations}} \footnote{Code tested using Mathematica 12.2.0.0 for Mac OS X} 
\begin{lstlisting}[breaklines]
z[t_] := t/Cosh[t]^2
K[s_] := 2^(s-1)/(1-2^(1-s))/Gamma[s+1]
Kp[s_] := Evaluate[D[K[s],s]]
f[s_] := Zeta[s]
Csgn[x_] := Sign[Re[x]]
Zetap[s_] := Evaluate[D[Zeta[s],s]]

alpha={0.048,0.235,0.8523,2.737}
alpha={48/1000,235/1000,8523/10000,2737/1000}
c = {0.0169,0.139,0.627,2.241}
c = {169/10000,139/1000,627/1000,2241/1000}
z0=57/100+157/100*I
R=1/10
inf=\[Infinity]
nj=4
n=1
Z1[s_]:=NIntegrate[z[t]*t^(s-1),{t,0,inf},WorkingPrecision->50, AccuracyGoal->5]
Z2[s_] :=NIntegrate[z[t/u1]*z[u1]/u1*t^(s-1),{u1,0,inf},{t,0,inf},WorkingPrecision->50, AccuracyGoal->5]
ZpZ0[s_] := NIntegrate[Log[t]*z[t]*t^(s-1),{t,0,inf},WorkingPrecision->50, AccuracyGoal->5]
ZpZ1[s_] := NIntegrate[Log[u1]*z[t/u1]*z[u1]/u1*t^(s-1),{u1,0,inf},{t,0,inf},WorkingPrecision->50, AccuracyGoal->5]

ExpApprox[x_] := Sum[1/Factorial[k]*x^k,{k,0,1}]
InvApprox[x_] := Sum[alpha[[j]]*Csgn[x]*Exp[-c[[j]]*x*Csgn[x]],{j,1,nj}] 
InvApprox2[x_] := Sum[alpha[[j]]*Csgn[x]*ExpApprox[-c[[j]]*x*Csgn[x]],{j,1,nj}] 
integrand1[\[Phi]_,k_]:=Sum[alpha[[j]]*Csgn[f[z0+R*Exp[I*\[Phi]]]]^(k+1)*(-1)^k/Factorial[k]*c[[j]]^k*Kp[z0+R*Exp[I*\[Phi]]]*K[z0+R*Exp[I*\[Phi]]]^k*Z1[z0+R*Exp[I*\[Phi]]]*I*R*Exp[I*\[Phi]]/(2*\[Pi]*I),{j,1,nj}] 
integrand2[\[Phi]_,k_]:=Sum[alpha[[j]]*Csgn[f[z0+R*Exp[I*\[Phi]]]]^(k+1)*(-1)^k/Factorial[k]*c[[j]]^k*Kp[z0+R*Exp[I*\[Phi]]]*K[z0+R*Exp[I*\[Phi]]]^k*Z2[z0+R*Exp[I*\[Phi]]]*I*R*Exp[I*\[Phi]]/(2*\[Pi]*I),{j,1,nj}] 
integrand3[\[Phi]_,k_]:=Sum[alpha[[j]]*Csgn[f[z0+R*Exp[I*\[Phi]]]]^(k+1)*(-1)^k/Factorial[k]*c[[j]]^k*K[z0+R*Exp[I*\[Phi]]]^(k+1)*ZpZ0[z0+R*Exp[I*\[Phi]]]*I*R*Exp[I*\[Phi]]/(2*\[Pi]*I),{j,1,nj}] 
integrand4[\[Phi]_,k_]:=Sum[alpha[[j]]*Csgn[f[z0+R*Exp[I*\[Phi]]]]^(k+1)*(-1)^k/Factorial[k]*c[[j]]^k*K[z0+R*Exp[I*\[Phi]]]^(k+1)*ZpZ1[z0+R*Exp[I*\[Phi]]]*I*R*Exp[I*\[Phi]]/(2*\[Pi]*I),{j,1,nj}] 
Table[N[I*R*Exp[I*\[Phi]]*Zetap[z0+R*Exp[I*\[Phi]]]/Zeta[z0+R*Exp[I*\[Phi]]]/(2*\[Pi]*I)] ,{\[Phi],0,2*\[Pi],2*\[Pi]/8}]
Table[N[I*R*Exp[I*\[Phi]]*Zetap[z0+R*Exp[I*\[Phi]]]*InvApprox[Zeta[z0+R*Exp[I*\[Phi]]]]/(2*\[Pi]*I)] ,{\[Phi],0,2*\[Pi],2*\[Pi]/8}]
Table[N[I*R*Exp[I*\[Phi]]*Zetap[z0+R*Exp[I*\[Phi]]]*InvApprox2[Zeta[z0+R*Exp[I*\[Phi]]]]/(2*\[Pi]*I)] ,{\[Phi],0,2*\[Pi],2*\[Pi]/8}]
Table[integrand1[\[Phi],0]+integrand2[\[Phi],1]+integrand3[\[Phi],0]+integrand4[\[Phi],1],{\[Phi],0,2*\[Pi],2*\[Pi]/8}]
\end{lstlisting}

\end{document}